\documentclass[12pt]{georeport}

\usepackage{amsmath}
\usepackage{amssymb}
\usepackage{amsthm}
\usepackage{amscd}
\usepackage{setspace}
\usepackage{wasysym}
\usepackage{graphicx}

\DeclareMathAlphabet{\mathcal}{OMS}{cmsy}{m}{n}
\newcommand{\bv}{\mathbf{v}}

\newcommand{\wl}{w_{\lambda}}
\newcommand{\dl}{d_{\lambda}}

\newcommand{\bR}{\mathbf{R}}

\newtheorem{theorem}{Theorem}

\setfigdir{Fig}

\begin{document}
\title{Wavefield Reconstruction Inversion: an example}
\author{William. W. Symes \thanks{The Rice Inversion Project,
Department of Computational and Applied Mathematics, Rice University,
Houston TX 77251-1892 USA, email {\tt symes@caam.rice.edu}.}}

\lefthead{Symes}

\righthead{1D WRI}

\begin{abstract}
  Nonlinear least squares data-fitting driven by physical process
  simulation is a classic and widely successful technique for the
  solution of inverse problems in science and engineering. Known as
  ``Full Waveform Inversion'' in application to seismology, it can
  extract detailed maps of earth structure from near-surface seismic
  observations, but also suffers from a defect not always encountered
  in other applications: the least squares error function at the heart
  of this method tends to develop a high degree of nonconvexity, so
  that local optimization methods (the only numerical methods feasible
  for field-scale problems) may fail to produce geophysically useful
  final estimates of earth structure, unless provided with initial
  estimates of a quality not always available. A number of alternative
  optimization principles have been advanced that promise some degree
  of release from the multimodality of Full Waveform Inversion,
  amongst them Wavefield Reconstruction Inversion, the focus of this
  paper. Applied to a simple 1D acoustic transmission problem, both
  Full Waveform and Wavefield Reconstruction Inversion methods reduce
  to minimization of explicitly computable functions, in an asymptotic
  sense. The analysis presented here shows explicitly how multiple
  local minima arise in Full Waveform Inversion, and that Wavefield
  Reconstruction Inversion can be vulnerable to the same
  ``cycle-skipping'' failure mode.
\end{abstract}

\section{Introduction}

Full waveform inversion (FWI) is the current nomenclature in the
seismology literature for data-fitting earth structure estimation
driven by wavefield modeling. The earth properties to be estimated
(material densities, stiffnesses, attenuation rates,...) form a vector
$c$ of spatially varying fields that appear as coefficients in systems of hyperbolic
partial differential equations, modeling seismic wave
propagation. The wavefields used in structure estimation (``imaging'')
are small motion disturbances of the earth's equilibrium state, so the
equations of motion are typically 
linear(ized). Right-hand side vectors $f$ in these systems model energy
input that initiates waves (earthquakes, man-made sources such as
explosives or mechanical vibrators). Data vectors $d$ are simulated by
sampling the solution fields at the locations of measurement devices
(accelerometers, microphones,...)  over appropriate time
intervals. The relation between the energy source $f$ and the
simulated data is linear in $f$, but nonlinear in the coefficient
vector $c$, so is naturally represented by a family of linear
operators $S[c]$ parametrized by the coefficient vector $c$.

The objective of FWI is to find $c$ and $f$, given $d$, so that
$S[c]f \approx d$. A typical method for achieving this goal is
the minimization of an objective misfit measure (objective, for
short), the most common choice being the square norm of a Hilbert
space in which the data is presumed to reside:
\begin{equation}
  \label{eqn:deffwi}
  \mbox{Given } d, \mbox{ find }c \mbox{ and }f \mbox{ to
    minimize }
  \|d -  S[c]f\|^2
\end{equation}
in which $\|\cdot\|$ is the norm in a suitable Hilbert space. This
approach was first suggested in the 1980's
(\cite[]{BamChavLai:79,Tara:84a,KolbColLai:86,Crasetal:90}, and many
other papers since then). Usually some form of regularization is
applied, to compensate for poorly determined aspects of $c$ and/or
$f$, as is explained in Tarantola's influential book
\cite[]{Tarantola:05}. Also, $f$ may be constrained in one way or
another to embody characteristics of field energy sources, or even
regarded as known (an example of this is given below).

Within a few years of its introduction into quantitative seismology,
FWI was understood to suffer from a severe limitation. Because of the
typical dimensions of earth models, and consequent cost of accurate
computation of $S$, iterative local optimization provides the only
feasible route to estimation of $c$ via solution of the optimization
problem \ref{eqn:deffwi}. However the objective function of this
optimization problem (the mean-square residual appearing in display \ref{eqn:deffwi}) has many local minima in general, most
having nothing to do with a usable estimate of earth structure. [An
explicit example of this multi-modal behaviour appears below.]
Reliable estimation of $c$ via iterative local optimization requires
that the initial estimate predict the correct arrival time of waves,
as they appears in the data, within a wavelength (in fact,
conventionally a half wavelength) at dominant
frequencies \cite[]{GauTarVir:86,VirieuxOperto:09}.

Despite this severe constraint, FWI has shown enough promise as a tool
for both industrial and academic seismology that it is now a
mainstream research topic, and to some extent a commercial
product. Estimation of sufficiently accurate initial models via
non-FWI methods is common practice, though what ``sufficiently
accurate'' means may be difficult to discern
\cite[]{Plessix:10}. The wavelength criterion mentioned in the last
paragraph may be made less onerous by collection of relatively
low-frequency data \cite[]{Wolfspar:16}. Finally, many alternatives to
straightforward least-squares data fitting have been suggested, some
of which appear to exhibit less tendency to develop local minima than
does the problem described in \ref{eqn:deffwi} \cite[]{Symes:09}.

This topic of this paper is one of these alternative approaches,
Wavefield Reconstruction Inversion (``WRI''). WRI was introduced by \cite{LeeuwenHerrmannWRI:13}, and further
developed by \cite{LeeuwenHerrmann:16,WangYingst:SEG16} and other
authors. It is based on the presumption that the correct source (right
hand side in the equations of motion) $q$ is
known, and combines a penalty for data misfit with a penalty for
failing to solve the equations of motion with the correct right hand side:
\begin{equation}
  \label{eqn:defwri}
  \mbox{Given } d \mbox{ and }q, \mbox{ find }c \mbox{ and }f \mbox{ to
    minimize }
  \|d -  S[c]f\|^2+\alpha^2\|f-q\|^2.
\end{equation}
\cite{Leeuwen2019note} points out that the formulation
\ref{eqn:defwri} is available not just for seismic problems, but for
any inverse problem based on a separable (partly linear) modeling
operator.

Examples given by \cite{LeeuwenHerrmannWRI:13} and others suggest that
the problem defined in display \ref{eqn:defwri} is less likely to
develop uninformative local minima than is the least squares problem
\ref{eqn:deffwi} in application to seismic inversion, if $\alpha$ is
sufficiently small. This might be so for several reasons: one is that
when less weight is put on making the residual in the wave equation
($f-q$) small, $f$ may be chosen to make the data residual $d-S[c]f$
small instead. One aspect of failure to predict the arrival times of
waves accurately is that small data residual is then difficult to
achieve. By maintaining fidelity to the data, local minimization of
the WRI objective might provide a spurious-minimum-free path to a
satisfactory estimate of $c$.

The main result of this paper is that in one simple case, in which all
of the necessary computations can be carried out by hand, this hope is
not realized: minimization of the WRI objective is just as likely to
be trapped in a spurious local minimizer as is the FWI objective.  The
context of this conclusion is a simple transmission inverse problem
for the 1D acoustic wave system, which models pulse transmission along
a 1D continuum from a source point to a receiver point. The pulses
used in this thought experiment are short, so the main information
content of the data is the time of transit from source to
receiver. The predominant information about the material model, in
this case reduced to the wave velocity (a scalar function of
position), is just this transit time, so I constrain both the target
wave velocity $c_*$ generating the data and the trial wave velocity
$c$ to be constant, that is, independent of position along the 1D
continuum. I introduce a family of inverse problems, depending on a
parameter $\lambda$ playing the role of wavelength. For sufficiently small
$\lambda$, it is possible to show explicitly that the FWI problem
possesses local minimizers far from the global minimizer at the target
$c_*$, and that initiating a local iterative optimization, such as
steepest descent or Newton's method, at a distance from $c_*$ bounded
below by a multiple of the ``wavelength'' $\lambda$ will result in
convergence to these spurious local minima. That is, FWI 
behaves in exactly the manner described in much of the literature on
this topic. However, an analysis of WRI applied to the same context
yields the same result: spurious local minima exist for sufficiently
small $\lambda$, and will be found by local optimization unless the
starting point is within $O(\lambda)$ of the target. That is, WRI
behaves in a manner qualitatively indistinguishable from FWI. In
particular, its ability to allow good fit to data for small $\alpha$
does not safeguard it from failure to converge globally to a ``good''
local minimum. In fact, the $\alpha \rightarrow 0$ limit of the WRI
objective is well-defined (after scaling by $1/\alpha^2$) and also behaves in
qualitatively the same way as the FWI objective. So the apparent
ability to maintain better data fit via reduction of $\alpha$ does not
lead to global behaviour asymptotically more amenable to local optimization.

This paper begins with a description of the 1D inverse transmission
problem and explicit computation of various components of the FWI
approach, based on explicit solution of the 1D acoustic system
presented in Appendix A. In the third section I introduce the $\lambda-$dependent
family of problems, and establish the asymptotic properties of FWI as
$\lambda \rightarrow 0$. The fourth section develops the algebraic
structure of WRI, culminating in a remarkable identity revealing WRI
to be equivalent to minimization of a weighted norm of the data
residual, with a weight operator depending on the coefficient vector
$c$. This identity has also been derived by \cite{Leeuwen2019note},
using a different argument. The result is quite general, applying to essentially any
realization of WRI, and already shows that its behaviour must be
closely related to that of FWI. In the fifth section I return to 
1D acoustics problem and the $\lambda$-dependent family of inverse
problems, compute that weight operator explicitly, and deduce the
global behaviour of WRI in this instance. The paper ends with a
discussion of the relation of the analysis presented here to a
previous analysis of optimization formulations of wave inversion
problems based on parameter-dependent quadratic forms
\cite[]{StolkSymes:03}, and a brief discussion of some other
alternatives to FWI.

\section{FWI for 1D Acoustics}
The example of FWI to be explored in this paper is one of the simplest possible, based on the 1D acoustics system connecting excess pressure $p$, particle velocity $v$, constitutive law defect (``source'') $f$, density $\rho$, and wave velocity $c$:
\begin{eqnarray}
\label{eqn:awe1d}
\frac{\partial p}{\partial t} + \rho c^2\frac{\partial 
  v}{\partial z} &=& f \nonumber\\
\rho \frac{\partial v}{\partial t} + \frac{\partial p}{\partial 
  z}&=&0\nonumber\\
 p,v&=&0, t \ll 0. 
\end{eqnarray}
The fields $p,v,f$ are functions of spatial position $z \in \bR$ and time $t \in \bR$, whereas $c, \rho$ are functions of $z$ alone, so that the system \ref{eqn:awe1d} is autonomous.

The system \ref{eqn:awe1d} has classical (smooth) solutions $(p,v)$
when $c, \rho,$ and $f$ are smooth, and $\log c$ and $\log \rho$ are
bounded on $\bR$, as is well-established \cite[]{Lax:PDENotes}. In
this paper, for reasons to be discussed below, $c$ and $\rho$ are
constrained to be constant ($z$-independent) in which case solutions
may be constructed by elementary methods (Appendix A). 

In fact, I shall assume $\rho>0$ to be fixed for the remainder of this
paper, that is, not updated in the inversion process. The wave
velocity $c$ will range over an interval: it is the parameter to be inverted.
To be specific, choose $c_{\rm max} > c_{\rm min} > 0$, and require
that $c$ satisfy $c_{\rm min} \le c \le c_{\rm max}$.

Limit observations to the time interval $[0,T]$, at the spatial (``receiver'')  location $z_r$. The modeling operator outputs the pressure trace $p(z_r,t)$ over the time interval $[0,T]$:
\begin{equation}
  \label{eqn:defmod}
  S[c]f = p|_{\{z_r\}\times [0,T]}
\end{equation}

The formulation of the inverse problem via least-squares requires a choice of Hilbert space structure for the domain and range of $S[c]$. A natural choice is
\begin{equation}
  \label{eqn:defmoddom}
  S[c]: L^2([z_{\rm min},z_{\rm max}] \times \bR) \rightarrow L^2[0,T]
\end{equation}
That is, the support of f will be assumed to lie in the strip $[z_{\rm min},z_{\rm max}] \times \bR$.
As it will turn out, the choice of spatial interval $[z_{\rm min},z_{\rm max}]$ is arbitrary, so long as it has positive length.

For homogeneous ($z$-independent) $c$, Appendix A provides an explicit expression:
\begin{equation}
  \label{eqn:defmodhom}
  S[c]f(t) = \frac{1}{2c}\int_{z_{\rm min}}^{z_{\rm max}} dz f\left(z,t - \frac{|z_r-z|}{c}\right) 
\end{equation}
From this expression it is simple to see that $S[c]$ is a bounded
operator with the domain and range described in
\ref{eqn:defmoddom}. The role of compact support in $z$ in ensuring
boundedness is also evident.

With this framework, the data is presumed to lie in the range space of
$S[c]$: $d \in L^2[0,T]$. It is simple to verify from the identity
\ref{eqn:defmodhom} that $S[c]$ is also surjective. That is, any data
at all can be fit by $S[c]$ with an appropriate choice of input. This
observation is important in the development of WRI, to be explained
below.

The version of FWI discussed here presumes that the source field
corresponding to the data $d$ is supported at a point $z_s \in \bR$,
and is known. The source position $z_s$ must satisfy $z_{\rm min} \le
z_s \le z_{\rm max}$ and $z_s \ne z_r$, but is otherwise arbitrary.

This point source field depends on a function of time (``wavelet'')
$w \in L^2(\bR)$. Formally, the resulting acoustic field satisfies the
system \ref{eqn:awe1d} with $f(z,t)=w(t)\delta(z-z_s)$. Inserting this
expression in the explicit expression \ref{eqn:defmodhom}, obtain
\begin{eqnarray}
  f(z,t) & = & w(t)\delta(z-z_s)\nonumber\\
  (S[c]f)(t) &=& \frac{1}{2c}w\left(t -
                              \frac{|z_r-z_s|}{c}\right)\nonumber\\
  &= &(S_p[c]w)(t)
  \label{eqn:defmodpt}
\end{eqnarray}
Thus for this restricted class of source fields, the FWI objective can be redefined as
\begin{equation}
  \label{eqn:deffwipt}
  J_{\rm FWI}[c,d.w] = \frac{1}{2}\|S_p[c]w -d\|^2.
\end{equation}
$S_p[c]$ is a bounded operator with domain $L^2(\bR)$ and range $L^2([0,T])$. 

To end this section, it is necessary to address an irritating
technical point: the point source defined above is not a member of the
domain of $S[c]$, as it was defined in display \ref{eqn:defmoddom}, so
the left-hand side of equation \ref{eqn:defmodpt} does not actually
make sense. The fix for this incompatibility actually elucidates the
relation between $S[c]$ and $S_p[c]$. Appendix B describes the
construction of a family of bounded injective operators  $E[c]: L^2(\bR) \rightarrow L^2([z_{\rm min},z_{\rm max}] \times \bR)$ for which
\begin{equation}
  \label{eqn:ext}
  S_p[c] = S[c] \circ E[c]
\end{equation}
This relation exhibits $S[c]$ as an {\em extension} of $S_p[c]$ as
described by \cite{Symes:09}. WRI for 1D acoustic transmission is
based on $S[c]$, and so is identified as an extended version of FWI.

The construction described in Appendix B requires that $z_s \in
(z_{\rm min},z_{\rm max})$, as mentioned above.

\section{Global Asymptotics of 1D transmission FWI}
The FWI objective is well-known to exhibit non-convexity unless
data frequency content is limited to a small range near 0 Hz, how
small being determined by other scales and by the extent to which the
initial wave velocity differs from the target.

To understand the non-convexity phenomenon and the
relation of the various scales, it is advantageous to introduce a
family of source wavelets, depending on a parameter $\lambda$, having
dimensions of time and playing
the role of wavelength:
\begin{equation}
  \label{eqn:wfam}
  \wl(t) = \frac{1}{\sqrt{\lambda}}w_1\left(\frac{t}{\lambda}\right).
\end{equation}
In the definition \ref{eqn:wfam}, the ``mother wavelet'' $w_1\in
C_0^{\infty}(0,1)$ has dimensionless argument, and the scaling is chosen so that
$\|\wl\|_{L^2(\bR)}$ is independent of $\lambda>0$.

Evidently there is no control of $c$ at all in the data if the time
interval of the observation, namely $[0,T]$, is so short that no
signal arrives within it. Accordingly, add to the other assumptions
made so far
the requirement that the transit time between source and receiver at
the slowest permitted velocity is less than $T$:
\begin{equation}
   \label{eqn:time}
  \frac{|z_s-z_r|}{c_{\rm min}} < T
\end{equation}

Choose a target wave velocity $c_* \in [c_{\rm min},c_{\rm
  max}]$. Introduce a family of consistent data $\dl$, generated by
$c_*$ and the wavelet family $\wl$:
\begin{equation}
  \label{eqn:dfam}
  \dl = S_p[c_*]\wl
\end{equation}
and a corresponding family of FWI objectives:
\begin{eqnarray}
  J_{\rm FWI}[c,\dl,\wl] &=& \frac{1}{2}\|\dl-S_p[c]\wl\|^2\nonumber\\
  &=&  \int_0^T\,dt\,\left|\frac{1}{2c_*}\wl\left(t-\frac{|z_s-z_r|}{c_*}\right)
    - \frac{1}{2c}\wl\left(t-\frac{|z_s-z_r|}{c}\right)\right|^2
  \label{eqn:fwiexpl}
\end{eqnarray}
Note that $\mbox{supp }\wl \subset [0,\lambda]$, so
\begin{equation}
  \label{eqn:wlsupps}
  \mbox{supp }S_p[c]\wl \subset \left[\frac{|z_s-z_r|}{c},\lambda
    +\frac{|z_s-z_r|}{c}\right] \cap [0,T]
\end{equation}
The transit time condition \ref{eqn:time} implies that there exists
$\lambda_0>0$ so that for $\lambda < \lambda_0$,
\[
  \lambda +\frac{|z_s-z_r|}{c} < T
\]
for all admissible $c$. That is, 
for $\lambda <\lambda_0$, $\mbox{supp }S_p[c]\wl \subset(0,T)\mbox{
  for }c \in [c_{\rm min},c_{\rm max}]$, and 
\begin{equation}
  \label{eqn:constnorm}
  \|S_p[c]\wl\|^2 = \frac{1}{4c^2}\int_{-\infty}^{\infty}\,dt\,\left|
    \wl\left(t-\frac{|z_s-z_r|}{c}\right)\right|^2 =
  \frac{\|w_1\|^2}{4c^2}.
\end{equation}

Recall that the object of this study is the global behaviour of
objective functions for velocity estimation: in this context, that
means the behaviour for $c$ far from $c_*$. Define
\begin{equation}
  \label{eqn:defL}
  L = \frac{2 c_{\rm max}^2}{|z_s-z_r|}.
\end{equation}
Then if $|c_*-c| > L\lambda$,
\begin{eqnarray}
  \left|\frac{|z_s-z_r|}{c}-\frac{|z_s-z_r|}{c_*}\right|
  &=&\frac{|c-c_*||z_s-z_r|}{cc_*} \nonumber\\
  &\ge& \frac{|c-c_*||z_s-z_r|}{c_{\rm max}^2}\nonumber\\
  &>& L\lambda\frac{|z_s-z_r|}{c_{\rm max}^2} \nonumber\\
  &=& 2\lambda.
      \label{eqn:shift}
\end{eqnarray}
That is, equation \ref{eqn:shift} shows that when the condition
\ref{eqn:defL} is satisfied, the infima of the supports of $S_p[c]\wl, \dl=S_p[c_*]\wl$
are further apart than the lengths of these supports. Then
necessarily
  $\mbox{supp }S_p[c]\wl \cap \mbox{supp }S_p[c_*]\wl = \emptyset$, so
$S_p[c]\wl$ and $S_p[c_*]\wl$ are orthogonal in $L^2[0,T]$,
and
\begin{equation}
  \label{eqn:fwiscrewed}
   J_{\rm FWI}[c,\dl,\wl]  =
  \frac{1}{2}\left(\frac{1}{4c^2}+\frac{1}{4c_*^2}\right)\|w_1\|^2.
\end{equation}
Amongst other consequences, one immediately deduces from the
expression \ref{eqn:fwiscrewed} the non-convexity result:

\begin{theorem}
  \label{thm:thm1}
  For $L>0$ given by equation \ref{eqn:defL} and $\lambda <
  \lambda_0$, the minimizer of $J_{\rm FWI}[c,\dl,\wl]$ on the
  complement of $[c_*-L\lambda, c_*+L\lambda]$ is $c=c_{\rm max}$.
\end{theorem}

That is, outside of a neighborhood of width proportional to a
wavelength, minimization of $J_{\rm FWI}$ yields a local minimizer far from
the target velocity $c_*$ that generates the (noise-free) data.

For this 1D problem, a happy 1D accident occurs: a descent
minimization starting at $c_0 < c_*$ will at least initially proceed
in the right direction. With sufficiently small steps, it is possible
that an interation might land in the (small) domain of attraction
around $c_*$. However neither this nor various other accidental
advantages stemming from the very special form of this problem should
be regarded as of any importance.
 
\section{Wavefield Reconstruction Inversion}
This section will describe Wavefield Reconstruction Inversion (WRI)
and develop some of its formal algebraic properties. The notation
$S[c]$ will represent a modeling operator based on wave dynamics of
some sort, depending on a vector $c$ of material parameters. The
conclusions developed here in fact apply to WRI in any such
setting. These conclusions will be applied to 1D acoustics in the
following section, with $c$ specialized to a scalar (wave velocity).

In the application below, the target source $q$ will be a point source
as in the previous section, however for the development of the basic
properties of WRI, that is immaterial.

Note that in general $\alpha$ is a dimensional parameter, having the
same dimensions as $S$.

\cite{LeeuwenHerrmannWRI:13} posed this problem slightly differently:
instead of the first order acoustics system (equation \ref{eqn:awe1d}
for the 1D case), they pose the wave dynamics in terms of the second
order wave equation for the pressure wavefield $p$. From this
viewpoint, $\partial(f-q) /\partial t$ is the residual, that is, the
difference between the image $\partial f/\partial t$ of the 2nd order
wave operator on $p$, and the assumed right-hand side
$\partial q/\partial t$. So in this form, the second term penalizes
the failure of $p$ to solve the wave equation with the assumed
source. The formulation presented here is equivalent, and was
introduced by \cite{WangYingst:SEG16}.

\cite{LeeuwenHerrmann:16} used the variable projection method
\cite[]{GolubPereyra:03}, eliminating the source $f$ in the inner step
and updating the bulk modulus (or an equivalent quantity) in the outer
step. That is, the problem \ref{eqn:defwri} is equivalent to
minimization of
\[
  J_{\rm WRI}^{\alpha}[c,d,q] =
  \mbox{min}_f \frac{1}{2}(\| d-S[c]f\|^2+\alpha^2\|f-q\|^2)
\]
over $c$. This is the approach that I shall pursue here.

It will turn out to be convenient to define the residual with the
target source $q$ as $r[c]=d-S[c]q$, and set $g=f-q$. Then this
definition can be rewritten as
\begin{equation}
\label{eqn:defvpm}
  J_{\rm WRI}^{\alpha}[c,d,q] =\mbox{min}_g \frac{1}{2}(\| r[c]-S[c]g\|^2+\alpha^2\|g\|^2)
\end{equation}

In any penalty method, control of the penalty
parameter has a large influence on the speed of
convergence. \cite{Aghamiry:19} use an augmented Lagrangian algorithm
to minimize the influence of the penalty weight choice. Alternatively, one can use a version of the discrepancy principle to adjust $\alpha$ dynamically \cite[]{FuSymes2017discrepancy}, as the WRI problem has the necessary features described in that paper.

``Most'' source fields $f$ are non-radiating, that is, $S[c]f=0$, and
such sources contribute nothing to the data fit term in the definition
of $J_{\rm WRI}^{\alpha}$. If the domain and range of $S$ were finite
dimensional (which of course they are, after discretization), then the
Fundamental Theorem of Linear Algebra identifies the null space of
$S[c]$ (the non-radiating sources) as the orthocomplement of the range
of the transpose $S[c]^T$ \cite[]{Strang:93}. In the
infinite-dimensional setting of this paper, the Closed Range Theorem
\cite[]{Yosida} states that the same is true if the range of $S[c]$ is
closed. There are various ways to ensure this property, but the
simplest is relevant here:
\begin{quote}
  For all admissible models $c$, $S[c]$ is surjective.
\end{quote}
That is, {\em any data can be fit exactly using the extended model
  space (the domain of $S$)}.
This property is a characteristic of of extended modeling
methods \cite[]{geoprosp:2008}: the ability to fit any data appears to
be a essential for such methods to produce objectives
without spurious local minima. It holds for the problems for which WRI
has been advocated. For the simple model problem considered
here, surjectivity follows from the explicit expression for $S[c]f$,
as was noted in the discussion following equation \ref{eqn:defmodhom}.

Assuming that $S[c]$ is surjective for any admissible $c$, the
orthocomplement of the subspace of non-radiating sources is the range
of the adjoint operator $S[c]^T$. In the
definition \ref{eqn:defvpm}, decompose $g = S[c]^Te + n$, in
which $e$ is the same type of object as $d$ and $S[c]n=0$ (that is,
$n$ is a non-radiating source), and note that the decomposition is
orthogonal. Then
\[
  J_{\rm WRI}^{\alpha}[c,d,q] =
  \mbox{min}_{e,n} \frac{1}{2}(\| d-S[c](S[c]^Te+q)\|^2+\alpha^2(\|S[c]^Te\|^2+\|n\|^2))
  \]
\begin{equation}
  \label{eqn:defvpmred}
 =  \mbox{min}_{e} \frac{1}{2}(\|r[c]-S[c]S[c]^Te\|^2+\alpha^2\|S[c]^Te\|^2).
 \end{equation}
 
This reformulation has some computational advantages \cite[]{WangYingst:SEG16,Herrmann:SEG19}, but also leads to a useful analytic transformation of the WRI problem. The minimizer on the RHS of equation \ref{eqn:defvpmred} is the solution $e=e_{\alpha}[c]$ of the normal equation
\[
 ( (S[c]S[c]^T)^2 + \alpha^2S[c](S[c]^T)e = S[c]S[c]^Tr[c]
\]
whence
\[
  S[c]S[c]^Te_{\alpha}[c] = S[c]S[c]^T(S[c]S[c]^T+\alpha^2I)^{-1}r[c]
\]
Since the null space of $S[c]$ is orthogonal to the range of $S[c]^T$
under the surjectivity assumption, $S[c]S[c]^T$ is injective, whence
\begin{equation}
  \label{eqn:norsol}
  e_{\alpha}[c]=(S[c]S[c]^T+\alpha^2I)^{-1}r[c]
\end{equation}
Consequently
\begin{eqnarray}
  J_{\rm WRI}^{\alpha}[c,d,q] &=&
  \frac{1}{2}(\|r[c]-S[c]S[c]^Te_{\alpha}[c]\|^2+\alpha^2\|S[c]^Te_{\alpha}[c]\|^2)\nonumber\\
  &=& \frac{1}{2}\left(\|r[c]-S[c]S[c]^T
    (S[c]S[c]^T+\alpha^2I)^{-1}r[c]\|^2\right. \nonumber\\
  & &
  \left.+\langle (S[c]S[c]^T+\alpha^2I)^{-1}r[c],
  S[c]S[c]^T (S[c]S[c]^T+\alpha^2I)^{-1}r[c]
  \rangle\right)  \nonumber\\
  &=&\frac{1}{2}\left(\|\alpha^2 (S[c]S[c]^T+\alpha^2I)^{-1}r[c]\|^2\right. \nonumber\\
  & &
  + \left. \alpha^2 \langle (S[c]S[c]^T+\alpha^2I)^{-1}r[c],
  - \alpha^2 (S[c]S[c]^T+\alpha^2I)^{-1}r[c]\rangle\right. \nonumber\\
  & &  + \left. \alpha^2\langle (S[c]S[c]^T+\alpha^2I)^{-1}r[c],r[c]\rangle\right) \nonumber\\
  &=&  \frac{\alpha^2}{2}\langle
      (S[c]S[c]^T+\alpha^2I)^{-1}r[c],r[c]\rangle.
   \label{eqn:longone}
\end{eqnarray}
Rearranging the RHS of equation \ref{eqn:longone}, obtain
\begin{equation}
  \label{eqn:defwrialt}
  J_{\rm WRI}^{\alpha}[c,d,q] = \frac{1}{2}\langle r[c], W_{\alpha}[c] r[c]\rangle
\end{equation}
with
\begin{equation}
  \label{eqn:defwriwt}
  W_{\alpha}[c] = \frac{\alpha^2}{2}(S[c]S[c]^T+\alpha^2I)^{-1}
\end{equation}
This remarkable identity shows that the WRI objective function is a
{\em weighted norm of the data residual $r[c]$}.

\cite{Leeuwen2019note} gives a different derivation of an identity
equivalent to equations \ref{eqn:defwrialt}, \ref{eqn:defwriwt}.

\section{Global Asymptotics of 1D transmission WRI}
The preceding section provides the necessary ingredients for an assessment of the relation between WRI and FWI. While the conclusion reached below applies to many wave propagation settings, the 1D acoustic setting is particularly simple and yet illustrates clearly the nature of this relation.

The first task is to give an explicit expression for the operator $S[c]S[c]^T$ appearing repeatedly in the expression \ref{eqn:defwriwt}. From the definition \ref{eqn:defmodhom}, it follows immediatlely that 
\begin{equation}
  \label{eqn:modhomtrans}
S[c]^Te(z,t)=
\left\{
    \begin{array}{c}
      \frac{1}{2c}e\left(t +  \frac{|z_r-z|}{c}\right), \,z_{\rm min} \le z \le z_{\rm max};\\
      0, else.
    \end{array}
  \right.
\end{equation}
whence
\[
  S[c]S[c]^Te(t) = \frac{z_{\rm max}-z_{\rm min}}{4c^2}e(t),
\]
that is,
\begin{equation}
  \label{eqn:modhomnormal}
  S[c]S[c]^T = \frac{z_{\rm max}-z_{\rm min}}{4c^2}I
\end{equation}
Thus the weight operator $W[c]$ appearing in \ref{eqn:defwrialt} takes
the form
\[
  W_{\alpha}[c] = u(c) I,
\]
\begin{equation}
  \label{eqn:wriwthom}
u(c) = \frac{\alpha^2}{2}\left(\frac{z_{\rm max}-z_{\rm min}}{4c^2}+\alpha^2\right)^{-1}.
\end{equation}

Next, suppose that $q = w\delta(z-z_s)$, that is, the target source is
a point source, so that $S[c]q = S_p[c]w$ in the notation used in
the discussion of FWI. Thus \ref{eqn:defwrialt} can be re-written as
\begin{equation}
  \label{eqn:wriscrewed}
  J_{\rm WRI}^{\alpha}[c,d,q]=u[c]J_{\rm FWI}[c,d,w]
\end{equation}

Recall the wavelength-dependent family of problems introduced in the derivation and
statement of Theorem \ref{thm:thm1}: target wave velocity $c_*$, wavelength parameter
$\lambda$, parametrized family of wavelets $\wl$ and corresponding
data $\dl$.

Define
\begin{equation}
  \label{eqn:defbeta}
\beta = \frac{z_{\rm max}-z_{\rm min}}{c_*^2} - 4 \alpha^2
\end{equation}
\begin{theorem}
  \label{thm:thm2}
  For $L$ as defined in \ref{eqn:defL}, and $\lambda <
  \lambda_0$, the minimizer of $J_{\rm WRI}^{\alpha}[c,\dl,\wl \delta(\cdot-z_s)]$ on the
  complement of $[c_*-L\lambda, c_*+L\lambda]$ is
  \begin{itemize}
  \item $c=c_{\rm max}$ if $\beta<0$;
  \item $c=c_{\rm min}$ if $\beta>0$;
  \item any $c < c_*-L\lambda$ or $>c_*+L\lambda$ if $\beta=0$.
  \end{itemize}
\end{theorem}
\begin{proof}
  From \ref{eqn:wriscrewed}, \ref{eqn:wriwthom}, and
  \ref{eqn:fwiscrewed}, if $|c-c_*| > L\lambda$,
\begin{eqnarray}
    J_{\rm WRI}^{\alpha}[c,\dl,\wl\delta(\cdot-z_s)] &=&
    \frac{\alpha^2}{2}\left(\frac{z_{\rm max}-z_{\rm
          min}}{4c^2}+\alpha^2\right)^{-1}\frac{1}{2}\left(\frac{1}{4c^2}+\frac{1}{4c_*^2}\right)\|w_1\|^2\nonumber\\
   & =&    \frac{\alpha^2}{4}\frac{1+\frac{c^2}{c_*^2}}{z_{\rm 
        max}-z_{\rm min}+4c^2\alpha^2}
  \label{eqn:yetanutta}
\end{eqnarray}
  The linear fractional function of $c^2$ on the RHS of equation
  \ref{eqn:yetanutta} is increasing, decreasing,
  or constant if $\beta>0, \beta<0$ or $\beta=0$, respectively.
\end{proof}

In other words, $J_{\rm WRI}^{\alpha}$ has local minima far from the target
velocity $c_*$, in the same way as does $J_{\rm FWI}$. One of the
local minima will be the result of a local optimization almost surely,
unless the initial estimate of $c$ is ``within a wavelength'' of the
target velocity.

Note that $L$ is independent of $\alpha$ (definition \ref{eqn:defL}),
and for small enough $\alpha$, $\beta > 0$ (definition
\ref{eqn:defbeta}). Conclude that the region $\{c \in [c_{\rm min},
c_{\rm max}]: |c-c_*| > L\lambda$ is
independent of $\alpha$, and the minimizer of
$J_{\rm WRI}[c,\dl,\wl\delta(\cdot-z_s)]$ in this region (away from
$c_*$) is $c=c_{\rm min}$ for small enough $\alpha$. Therefore taking
$\alpha$ small does not change the multimodal nature of
$J_{\rm WRI}^{\alpha}$: there remain multiple far-apart local minima,
no matter how small $\alpha$ may be.

\section{Discussion}
Theorems \ref{thm:thm1} and \ref{thm:thm2} call out the chief
conclusions of this work: that at least for the 1D acoustic
transmission inverse problem, both $J_{\rm FWI}[c,d,w]$ and
$J_{\rm WRI}^{\alpha}[c,d,q]$ exhibit local minima (in $c$) far from the global
minimum for consistent data, the domain of attraction of the
global minimizer can be arbitrarily small, and these properties
persist as $\alpha \rightarrow 0$. These are striking
conclusions, but the 1D acoustic transmission problem is very special
and lacks fidelity to field practice. I shall show how these
approaches to solving this special problem share properties with a
much larger family of inversion methods. The theory developed to
explain these properties suggests methods that may not suffer the
non-convexity of FWI and WRI, and in fact do not in several cases that
I will mention. 

First, a consequence of the results proven here: $J_{\rm FWI}[c,d,w]$
and $J_{\rm WRI}^{\alpha}[c,d,q]$ are not smooth as joint functions of model
($c$) and data ($d$) vectors. If they were, their derivatives would be
bounded uniformly over bounded sets in $c,d$, but the two main results
show that this is not the case. As $\lambda \rightarrow 0$, $d=\dl$
varies within a ball $B \subset L^2([0,T])$ of radius $\|d_1\|$ (since
the $\|\dl\|$ is independent of $\lambda$), but the value of either
objective changes from a positive value (bounded away from zero
independently of $\lambda$) to zero over an interval of $c$ of length
$O(\lambda)$. Therefore the derivatives with respect to $c$ of both
objective functions are not bounded over
$[c_{\rm min},c_{\rm max}] \times B$.

While lack of smoothness is not in itself the most important property
established in the preceding sections, it is a necessary condition for
stable and reliable parameter recovery via local
optimization. Moreover, necessary conditions for smoothness are
known for a much wider class of quadratic form objectives for inverse
problems. 
  
These results
concern optimization problems of the general form
\begin{equation}
  \label{eqn:gen1}
  \mbox{Given }d\mbox{, find }c\mbox{ to extremize }
   J[c,d] = \langle G[c]d, A[c] G[c]d \rangle
\end{equation}
In this prescription, $d$ is a data vector, as in the examples above,
$c$ is a vector of material parameters to be estimated, $G[c]$ is a
$c$-dependent family of operators, whose common range is a Hilbert
space with inner product $\langle \cdot, \cdot \rangle$.
$A[c]$ is an operator-valued function of $c$, with domain and the
range equal to the range of $G[c]$.

\cite{StolkSymes:03} assume that the
operator-valiued function $G[c]$ is of a class typical of modeling
operators for wave equation inverse problems, or their inverses or
adjoints. The precise characterization of these so-called {\em microlocally elliptic Fourier Integral
Operators} is quite technical \cite[]{Dui:95}. Roughly speaking, such operators map high-frequency localized wave packets to other
such packets with well-defined changes of position and direction of
oscillation. The simulation operators $S[c]$ and $S_p[c]$ figuring in the preceding
discussion are particularly simple examples of this type.

If the dependence of such an operator $G[c]$ on $c$ is of sufficiently full rank, in
the sense that destination packets can be shifted in any direction by
changing $c$ appropriately, along with a couple of other technical
assumptions, one can conclude that $J$ as defined in \ref{eqn:gen1} is
smooth in $c$ and $d$ jointly if and only if the operator $A[c]$ is a 
{\em pseudodifferential operator} - again, a class of operators whose
precise definition is quite technical \cite[]{Dui:95,Tay:81}. However these operators also
have a rough characterization: they do not change the location of
oscillatory wave packets
or alter their direction of oscillation, but only scale such packets
by smooth functions, to good approximation. 

With a bit of fiddling, the FWI problem \ref{eqn:deffwipt} for 1D 
acoustic transmission inversion can be rewritten in the form \ref{eqn:gen1}. Note
that $S_p[c]$ is invertible (more specifically, has a right inverse; if [0,T] were
extended to $(-\infty,\infty)$ it would have a left inverse too). From
the definition \ref{eqn:dfam} of the data family $\dl$, one sees that
$\wl = S_p[c_*]^{-1}\dl$. Therefore
\begin{eqnarray}
  \label{eqn:anutta}
  J_{\rm FWI}[c,\dl,\wl] &=&
                             \frac{1}{2}\|(I-S_p[c]S_p[c_*]^{-1})\dl\|^2 \nonumber\\
 & =&\frac{1}{2}(\|\dl\|^2 + \|S_p[c]S_p[c_*]^{-1}\dl\|^2) + \langle \dl,
  S_p[c]S_p[c_*]^{-1}\dl\rangle
\end{eqnarray}
For the operator family $S_p[c]$ defined above, it is easy to see that
the second term in the right hand side of equation \ref{eqn:anutta} is
smooth in $c$, and the first is constant. The third can be rewritten
as
\begin{equation}
  \label{eqn:evenmore}
  \langle \dl,  S_p[c]S_p[c_*]^{-1}\dl\rangle = \langle S[c]^T\dl, 
  (S_p[c]^T S_p[c_*])^{-1}S_p[c]^Td\rangle.
\end{equation}
The RHS of equation \ref{eqn:evenmore} has the form \ref{eqn:gen1} with the choices $G[c]=S_p[c]^T$,
$A[c]=(S_p[c]^T S_p[c_*])^{-1}$. A similar manipulation exhibits
$J_{\rm WRI}$ as the sum of harmless terms and a quadratic form of the
form \ref{eqn:gen1}.

Given the rough understanding of the results of \cite{StolkSymes:03}
sketched above, one would conclude that neither $J_{\rm FWI}$ nor
$J_{\rm WRI}$ are likely to be smooth jointly in $c$ and $d$. Indeed,
apart from scale, $S_p[c]$ is composition with a shift (translation)
by $(z_{\rm max}-z_{\rm min})c^{-1}$, so
$A[c]=(S_p[c]^T S_p[c_*])^{-1}$ is composition with a shift by
$(z_{\rm max}-z_{\rm min})(c^{-1}-c_*^{-1})$. Thus application of
$A[c]$ does not leave the position of a wave packet fixed, unless
$c=c_*$ - and indeed $A[c]$ is not a pseudodifferential operator
unless $c=c_*$. On the other hand, $G[c]$ is a shift operator, the
simplest prototype of an elliptic Fourier Integral Operator. Therefore
the conclusion, derived directly from Theorem \ref{thm:thm1}, that
$J_{\rm FWI}$ is not smooth jointly in $c$ and $d$ would also appear
to follow from the main result of \cite[]{StolkSymes:03}. This
conclusion can be made precise by proper attention to detail, and the
same is true of $J_{\rm WRI}$.

In the context of the 1D acoustic transmission problem as formulated
here, the question immediately arises: do quadratic forms \ref{eqn:gen1}
exist that are smooth jointly in $c$ and $d$, and whose global
minimizer is the correct velocity $c=c_*$? An affirmative answer is
provided for precisely this example problem in
\cite[]{wwsorcas:19-02}. The operator $G[c]=S_p[c]^T$ is precisely the
same as appeared in the reformulation of $J_{\rm FWI}$. Ignoring a
$c-$dependent multiplier,
\[
  A[c]u(t) = tu(t).
\]
Applied to the wavelength-dependent family of data $\dl$ and source
wavelets $\wl$ used repeatedly throughout this paper, $A[c]$ yields a
vanishing result as $\lambda \rightarrow 0$ for the correct velocity
$c=c_*$ and a stably non-zero result otherwise. For these choices, it
can be established that
\[
  \frac{d}{dc}J[c,\dl]
  \left\{\begin{array}{c}
           > 0 \mbox{ if }c<c_*+O(\lambda),\\
           < 0 \mbox{ if }c>c_*+O(\lambda)
         \end{array}
       \right.
     \]
That is, all local minima of $J[c,\dl]$ lie within $O(\lambda)$ (``a
wavelength'') of the target velocity $c_*$. Not surprisingly the
wavelength parameter also regulates the accuracy of the inversion.

The reader is directed to \cite[]{HuangNammourSymesDollizal:SEG19} for
an extensive discussion of other similar source extension methods for
various wave inversion problems, and for
references to earlier work on this topic.

\section{Conclusion}
The tendency of iterative FWI to become trapped in uninformative local
minima has been much discussed and still drives a substantial
worldwide research program, almost 35 years after the phenomenon was
first identified. WRI is amongst the many remedies proposed for this
pathology, and numerical experiments have appeared to suggest that it
may succeed. The example investigated in this report is simple enough
to allow for rigorous mathematical conclusions regarding the behaviour
of both FWI and WRI. The complete explanation for the behaviour of FWI
is no surprise. As it turns out, the same conclusion may be reached
for WRI: in this example at least, iterative minimization of the WRI
objective it is no more likely to produce a useful estimate of wave
velocity than is FWI, and for the same reason - indeed, the two are
very closely linked (equations \ref{eqn:defwrialt},
\ref{eqn:defwriwt}).

While these specific conclusions are of course tied to the extremely
simple homogeneous acoustic 1D transmission inverse problem studied
here, the relations \ref{eqn:defwrialt}, \ref{eqn:defwriwt} are
straightforward algebraic properties of WRI and appear to link it
closely to FWI in any wave propagation setting. As explained in the
discussion section, even mere smoothness of a quadratic form objective
function in both the model and data parameters may impose restrictions
on the operators involved in the construction of the form. These
restrictions are generally not met by any version of FWI. An
examination of other versions of WRI from this point of view may prove
informative.


\append{1D Radiation Problem}
Begin with the 1D acoustics point source system. 
\begin{eqnarray}
\label{eqn:awe1dptsrc}
\frac{\partial p}{\partial t} +\rho c^2\frac{\partial 
  v}{\partial z} &=& w(t)\delta(z-z_s) \nonumber\\
\rho \frac{\partial v}{\partial t} + \frac{\partial p}{\partial 
  z}&=&0\nonumber\\
 p,v&=&0, t \ll 0. 
\end{eqnarray}
Since the right hand side is singular, so is the solution, so it must
be a solution in the weak sense. It follows from the weak solution
conditions that the pressure is continuous at $z=z_s$, whence $v$ must
have a discontinuity. 

In $z \ne z_s$, the right hand side 
vanishes, so the solution must be locally a combination of plane
waves; causality implies that
\[
p(z,t)=a\left(t -\frac{|z-z_s|}{c}\right), \, v(z,t)=\mbox{sgn}(z-z_s) b\left(t -
  \frac{|z-z_s|}{c}\right)
\]
From the second dynamical equation (Newton's law) it follows that $b =
a/(\rho c)$. The singularity on the LHS of the first dynamical
equation (constitutive law) is
\[
\rho c^2 [v]_{z=z_s}\delta(z-z_s) =
2\rho c^2 b\delta(z-z_s) = 2c a\delta(z-z_s).
\] 
This must in turn equal the RHS of the constitutive law, whence
$a=w/(2c)$. Thus
\begin{eqnarray}
\label{eqn:sol1dptsrc}
p(z,t) &=& \frac{1}{2c}w\left(t - \frac{|z-z_s|}{c}\right) \nonumber \\
v(z,t) &=& \mbox{sgn}(z-z_s)\frac{1}{2\rho c^2}w\left(t -\frac{|z-z_s|}{c}\right)
           \nonumber \\
\end{eqnarray}
This result (computation of the Green's function for the acoustic
system) permits an explicit expression for the system with a
space-time source:
\begin{eqnarray}
\label{eqn:awe1d}
\frac{\partial p}{\partial t} +\rho c^2\frac{\partial 
  v}{\partial z} &=& f(z,t) \nonumber\\
\rho \frac{\partial v}{\partial t} + \frac{\partial p}{\partial 
  z}&=&0\nonumber\\
 p,v&=&0, t \ll 0. 
\end{eqnarray}
Since
\[
  f(z,t) = \int dz_1\,f(z_1,t)\delta(z-z_1)
\]
obtain
\begin{eqnarray}
\label{eqn:sol1dp}
p(z,t) &=& \frac{1}{2c}\int dz_1 f\left(z_1,t -
           \frac{|z-z_1|}{c}\right) \\
  \label{eqn:sol1dv}
v(z,t) &=& \frac{1}{2\rho c^2} \int dz_1 \mbox{sgn} (z-z_1) f\left(z_1,t - \frac{|z-z_1|}{c}\right)
\end{eqnarray}

\append{Equivalence of point and non-point sources}

As noted in the text, the point source $w(t)\delta(z-z_s)$ is not a
member of the domain of the simulation operator $S[c]$, as it is not
square-integrable. The object of this appendix is to construct a square-integrable right-hand side in the
system \ref{eqn:awe1d} for which the pressure field $p$
is the same as that of the weak solution to the point source problem \ref{eqn:awe1dptsrc}
constructed in the last section, near the receiver point $z=z_r$, and
to exhibit this square-integrable replacement for the point source as
the image of the point source wavelet under a bounded extension map,
as in equation \ref{eqn:ext}.

One step in this construction involves building a constitutive defect
(pressure) source that is equivalent to a force (velocity) source, in the sense of generating the
same solution outside of the source support. This construction is
presented here in the context of the 1D acoustic problem, but is a
special case of a much more general construction of considerable
interest in its own right \cite[]{BurridgeKnopoff:64}.

Let
$\epsilon$ be any positive number $<|z_r-z_s|$. 
Denote by $(p,\bv)$ the (weak) solution \ref{eqn:sol1dptsrc} of the point source
problem constructed in the last section.
Pick $\phi \in C_0^{\infty}(\bR)$ so that
$\phi = 1$ if $|z-z_s| \le \epsilon/2$ and $\phi(z)=0$ if $|z-z_s|
\ge\epsilon$. Set $p_0=p(1-\phi)$, $v_0=v(1-\phi)$. Then
\begin{eqnarray}
\label{eqn:awecut}
\frac{\partial p_0}{\partial t} + \rho c^2 \frac{\partial
  v_0}{\partial z} &=&
                            f_0  \nonumber\\
\rho \frac{\partial v_0}{\partial t} +\frac{\partial p_0}{\partial z} &=& g_0 \nonumber\\
\end{eqnarray}
in which
\begin{eqnarray}
  \label{eqn:rhscut}
  f_0(z,t) &=& -\rho c^2 v(z,t) \frac{\partial}{\partial z}(1-\phi(z))\nonumber\\
  g_0(z,t) &=& -p(z,t) \frac{\partial}{\partial z}(1-\phi(z))
\end{eqnarray}
vanish near $z=z_s$. If $w \in L^2(\bR)$ and vanishes for large
negative $t$ (as it must, for the system \ref{eqn:awe1dptsrc} to be
compatible), then from expressions \ref{eqn:sol1dptsrc} the
distributions $p,v$ are locally square-integrable in
$\{z:|z-z_s| \ge \epsilon/2\} \times \bR$ and vanish for large
negative $t$, whence the same is true of $f_0,g_0$.

Assume for the moment that $w \in C^{\infty}_0(\bR)$, so that $p,v$
are smooth away from $z=z_s$ and $p_0, v_0, f_0, g_0$ are smooth. Then
$p_0$ is also the solution of the second-order initial value problem
\begin{eqnarray}
  \label{eqn:awe2ord}
  \frac{1}{\rho c^2} \frac{\partial^2 p_0}{\partial t^2} -
  \frac{1}{\rho}\frac{\partial^2 p_0}{\partial z^2} &=& F\nonumber \\
  p &=& 0, t\ll 0
\end{eqnarray}
   
with the right-hand side $F$ given by
\begin{equation}
  \label{eqn:rhs2ord}
  F = \frac{1}{\rho c^2}\frac{\partial f_0}{\partial t} -
  \frac{1}{\rho}\frac{\partial g_0}{\partial z}
\end{equation}
Define $f$ by
\begin{equation}
  \label{eqn:defequivrhs}
  f(z,t) = \rho c^2 \int_{-\infty}^t \,ds\,F(z,t) =  f_0(z,t) - c^2\int_{-\infty}^{t}\,ds\,\frac{\partial
    g_0}{\partial z}(z,s)
\end{equation}
Then setting
\begin{eqnarray*}
  p_1 &=& p_0,\\
  v_1 &=& \frac{1}{\rho}\int_{-\infty}^t\frac{\partial p_0}{\partial z}
\end{eqnarray*}
it follows from \ref{eqn:awe2ord} and \ref{eqn:defequivrhs} that
$p_1,v_1$ solves \ref{eqn:awe1d} with $f$ as given above. Since
$p_1=p_0$, and $p_0=p$ in a neighborhood of $z=z_r$, it follows that
\begin{equation}
  \label{eqn:impext}
  S[c]f = S_p[w],
\end{equation}
that is, that using RHS $f$ in \ref{eqn:awe1d} produces the same
pressure field near $z=z_r$ as does the point source in
\ref{eqn:awe1dptsrc}. Also
\[
  f(z,t) = -\rho c^2 v(z,t) \frac{\partial}{\partial z}(1-\phi(z)) - c^2
  \int_{-\infty}^t\,ds\, \frac{\partial}{\partial z}\left(-p(z,s) \frac{\partial}{\partial z}(1-\phi(z))\right)
\]
\[
  = - c^2 \left(\int_{-\infty}^t\,ds\,\left( \rho \frac{\partial v}{\partial
      t} -\frac{\partial p}{\partial
      z}\right)(z,s)\frac{\partial}{\partial z}(1-\phi(z))
  -p(z,s)\frac{\partial^2}{\partial z^2}(1-\phi(z)) \right)
\]
\[
  = -2\rho c^2 v(z,t)\frac{\partial}{\partial z}(1-\phi(z)) +
  c^2\frac{\partial^2}{\partial z^2}(1-\phi(z))\int_{-\infty}^t
  \,ds\,p(z,s)
\]
using the second equation (momentum balance) in the system
\ref{eqn:awe1dptsrc}. Use \ref{eqn:sol1dptsrc} to replace $p,v$ by
explicit expressons in $w$:
\[
    = -\mbox{sgn}(z-z_s)w\left(t -\frac{|z-z_s|}{c}\right)\frac{\partial}{\partial z}(1-\phi(z)) +
  \frac{\partial^2}{\partial z^2}(1-\phi(z))\frac{c}{2}\left(\int_{-\infty}^tw\right)\left(t - \frac{|z-z_s|}{c}\right)
\]
\begin{equation}
  \label{eqn:defextop}
  =E[c]w(z,t)
\end{equation}
whence the image of $w$ under $E[c]$ is square-integrable,
$E[c]$ extends to a bounded operator $L^2(\bR) \rightarrow
L^2([z_{\rm min},z_{\rm max}]\times \bR)$, and from equation \ref{eqn:impext}
\begin{equation}
  \label{eqn:extapp}
S[c]\circ E[c] = S_p[c]
\end{equation}
as asserted in equation \ref{eqn:ext}.

\bibliographystyle{seg}
\bibliography{../../bib/masterref}\end{document}